\documentclass{siamart190516}

\usepackage{amsmath,amssymb,amsfonts,cmmib57,mathdots,xcolor,mathrsfs}
\usepackage[latin1]{inputenc}
\usepackage{enumerate}
\usepackage{caption}
\usepackage{subcaption}

\newcommand{\C}{\mathbb{C}}

\newcommand{\UU}{\mathcal{U}}
\newcommand{\VV}{\mathcal{V}}

\DeclareMathOperator{\coker}{coker}

\newtheorem{remark}[theorem]{Remark}
\newtheorem{example}[theorem]{Example}

\title{Perturbation theory of transfer function matrices}
\author{Vanni Noferini\thanks{Aalto University, Department of Mathematics and Systems Analysis, P.O. Box 11100, FI-00076, Aalto, Finland. (\texttt{vanni.noferini@aalto.fi})  Supported by an Academy of Finland grant (Suomen Akatemian p\"{a}\"{a}t\"{o}s 331230).}
\and
Lauri Nyman\thanks{Aalto University, Department of Mathematics and Systems Analysis, P.O. Box 11100, FI-00076, Aalto, Finland. (\texttt{lauri.s.nyman@aalto.fi}).}
\and
Javier P\'{e}rez\thanks{University of Montana, Department of Mathematical Science, 32 Campus Dr, Missuola, MT 59812, United States. (\texttt{javier.perez-alvaro@mso.umt.edu}) }
\and
Mar\'{i}a C. Quintana\thanks{Corresponding author. Aalto University, Department of Mathematics and Systems Analysis, P.O. Box 11100, FI-00076, Aalto, Finland. (\texttt{maria.quintanaponce@aalto.fi})  Supported by an Academy of Finland grant (Suomen Akatemian p\"{a}\"{a}t\"{o}s 331230) and by the Agencia Estatal de Investigaci\'{o}n of Spain through grant PID2019-106362GB-I00 MCIN/ AEI/10.13039/501100011033/.}}
\begin{document}
\maketitle
\slugger{simax}{xxxx}{xx}{x}{x--x}
\renewcommand{\thefootnote}{\fnsymbol{footnote}}

\newcommand{\la}{\lambda}

\renewcommand{\thefootnote}{\arabic{footnote}}

\begin{abstract} Zeros of rational transfer function matrices $R(\la)$ are the eigenvalues of associated polynomial system matrices $P(\la)$, under minimality conditions. In this paper we define a structured condition number for a simple eigenvalue $\la_0$ of a (locally) minimal polynomial system matrix $P(\la)$, which in turn is a simple zero $\la_0$ of its transfer function matrix $R(\la)$. Since any rational matrix can be written as the transfer function of a polynomial system matrix, our analysis yields a structured perturbation theory for simple zeros of rational matrices $R(\la)$. To capture all the zeros of $R(\la)$, regardless of whether they are poles or not, we consider the notion of root vectors. As corollaries of the main results, we pay particular attention to the special case of $\la_0$ being not a pole of $R(\la)$ since in this case the results get simpler and can be useful in practice. We also compare our structured condition number with Tisseur's unstructured condition number for eigenvalues of matrix polynomials, and show that the latter can be unboundedly larger. Finally, we corroborate our analysis by numerical experiments.
\end{abstract}

\pagestyle{myheadings}
\thispagestyle{plain}

\section{Introduction} Given a rational matrix $R(\la)$, the Rational Eigenvalue Problem (REP) is often defined as the problem of finding scalars $\la_0$ such that $R(\la_0)$ has finite entries (that is, $\la_0$ is not a pole of $R(\la)$) and that there exist nonzero constant vectors $x$ and $y$ (called eigenvectors) satisfying $R(\la_0)x=0$ and $y^{T}R(\la_0)=0,$ under the assumptions that $R(\la)$ is regular, i.e., $R(\lambda)$ is square and its determinant is not identically equal to zero.  However, zeros of rational matrices can also be poles: a situation not uncommon in certain applications such as control theory \cite{Kailath,vannof,Rosen70}. Thus, in this paper we will define the REP regardless of whether a zero is also a pole or not.  For that, we extend the definition of eigenvectors by using the more general notion of root vectors \cite{DN20,GT,vannof}. That is, we define the REP as the problem of finding scalars $\la_0$ and rational vectors $x(\la)\in\mathbb{C}(\lambda)^{m}$ and $y(\la)\in\mathbb{C}(\lambda)^{m}$ with $x(\la_0)\neq 0$ and $y(\la_0)\neq 0$ such that 
$\lim_{\lambda \to \lambda_0}R(\lambda)x(\la) = 0$ and $ \lim_{\lambda \to \lambda_0} y(\la)^T R(\lambda)=0.$ The above definition captures all the zeros of regular rational matrices $R(\la)$.

It is known \cite{Rosen70} that any rational matrix $R(\la)$ can be written (or is directly given) as the transfer function matrix of a polynomial system matrix. That is, of the form 
\begin{equation}\label{eq:general representation}
R(\lambda) = D(\la)+C(\la)A(\la)^{-1}B(\la),
\end{equation}
where $A(\la)$, $B(\la)$, $C(\la)$ and $D(\la)$ are arbitrary polynomial matrices \cite{GohLR09}, with $A(\la)$ regular. Clearly, the expression in \eqref{eq:general representation} is not uniquely determined by $R(\la)$. Then, under minimality conditions \cite{local,Rosen70}, the zeros of $R(\la)$ are the eigenvalues of the associated polynomial matrix
\begin{equation}\label{eq:psm}
P(\la)=\begin{bmatrix}
-A(\la) & B(\la)\\
C(\la) & D(\la)
\end{bmatrix},
\end{equation}
which is said to be a polynomial system matrix of $R(\la)$. In this paper, we will study the structured conditioning of simple zeros  of \eqref{eq:psm}, allowing perturbation of each block; if \eqref{eq:psm} is (locally) minimal, the analysis thus provides a structured perturbation theory of simple zeros of rational matrices written as in \eqref{eq:general representation}.

Rational matrices appear directly from applications, as in linear systems and control theory, or as approximations to nonlinear eigenvalue problems (NLEPs) \cite{nlep,automatic}. An example of rational matrix arising from application \cite{MV04} is the following:
\begin{equation}\label{eq:real_application}
R(\la)=-K+\lambda M + \la^{2}\sum_{i=1}^k \frac{1}{\omega_i- \la}C_i,
\end{equation}
which can be easily written as the transfer function of a polynomial system matrix, for example, as follows:
\begin{equation}\label{eq:real_application_representation}
R(\la)=-K+\lambda M + \begin{bmatrix}\la I \cdots \la I
\end{bmatrix}\begin{bmatrix}
(w_1-\la) I & & \\
& \ddots & \\
& & (w_k-\la) I
\end{bmatrix}^{-1}\begin{bmatrix}
\la C_1 \\
\vdots \\
\la C_k
\end{bmatrix}.
\end{equation}
Another example is the problem \texttt{loaded string} from \cite{NLEVP}:
\begin{equation}\label{eq:real_application2}
R(\la)=A-\la B + \frac{k\la}{\la-k/m} e_n e_n^T,
\end{equation}
where $A,B$ are certain $n\times n$ matrices, $k,m$ are parameters, and $e_n$ denotes the $n$th column of the $n\times n$ identity matrix. The rational matrix \eqref{eq:real_application2} is also easy to write as a transfer function, for instance:
\begin{equation}\label{eq:real_application_representation2}
R(\la)=A-\la B + \la e_n \begin{bmatrix}
\la-k/m
\end{bmatrix}^{-1} k e_n^{T}.
\end{equation}

\subsection{Condition numbers}\label{sub:conddef}

The condition number of a function measures how much the output of the function can change for a small perturbation in the input. To establish a general framework, consider two normed vector spaces $(\UU,\| \cdot \|_\UU)$ and $(\VV,\| \cdot \|_\VV)$  and a function $f : \UU \rightarrow V,\, u  \mapsto f(u) $. The (worst-case) \emph{absolute condition number of $f$ at $u$} is defined as \cite{Rice}
\begin{equation}\label{eq:abscond}
\kappa_f (u)= \lim_{\epsilon \rightarrow 0} \sup_{\| e \|_\UU \leq 1} \frac{\| f(u+\epsilon e)-f(u)\|_\VV}{\epsilon}.
\end{equation}
If $f$ is Fr\'{e}chet differentiable at $u$, then this definition implies that $\kappa$ is the operator norm of the Fr\'{e}chet derivative of $f$ at $u$. However, \eqref{eq:abscond} is valid more generally. The (worst-case) \emph{relative condition number of $f$ at $u$} is defined analogously as 
\begin{equation}\label{eq:relcond}
\kappa_{\mathrm{f,rel}} (u) =
\lim_{\epsilon \rightarrow 0} \sup_{\frac{\| e \|_\UU}{\|u\|_\UU} \leq 1} \frac{\| f(u+\epsilon e)-f(u)\|_\VV}{\epsilon \| f(u)\|_\VV}
= \kappa_f(u) \cdot \frac{\| u \|_\UU}{\| f(u) \|_\VV}.
\end{equation}
Condition numbers are geometric invariants, meaning that (once an arbitrary, but fixed, choice of norms is made) they are an intrinsic property of the function $f$ and the point $u$; in other words, they do not depend on the algorithm used in practice to compute $f(u)$. In numerical analysis, their use is popular to measure how difficult it is to accurately compute $f$ in finite precision arithmetic. This has been effectively summarized by N. Higham \cite[Sec. 1.6]{asna} with the ``rule of thumb" 
$$ \mathrm{forward} \ \mathrm{error} \lesssim \kappa \cdot \mathrm{backward} \ \mathrm{error},$$
which comes from a first order expansion and therefore is technically true only for finite condition numbers and small enough values of the backward error. In practice, it is often still valid more generally as an approximate inequality. For a backward stable algorithm, the backward error is guaranteed to be small, but whether this implies that the output of the computation is accurate (that is, whether the forward error is small) also depends on the condition number $\kappa$. Finally, \emph{structured} condition numbers can be defined analogously, but with the further restriction that the set of allowed perturbations $e$ must satisfy certain properties. In other words, instead of taking the supremum over the unit ball $\| e \|_\UU \leq 1$ as in \eqref{eq:abscond}, for defining a structured condition number one takes the supremum over a certain \emph{subset} of the unit ball. For example, if $u$ is a Hermitian matrix and $\UU = \C^{n \times n}$, then one may impose that $e$ is such that $u+\epsilon e$ is also a Hermitian matrix. Note that, by the properties of the supremum, this immediately implies that a structured condition number is less than or equal to the corresponding unstructured condition number; equality is however generally possible.

In this paper, we specialize the general theory to condition numbers of simple zeros (also called eigenvalues) of either rational or polynomial matrices, always assumed to be regular. Suppose that $\la_0$ is a simple zero for a rational matrix $R(\la)$, associated with left and right root vectors \cite{GT,vannof} $y(\la)$ and $x(\la)$ respectively. Then, as we prove in Lemma \ref{lem_der}, for any sufficiently small $\epsilon >0$ and any perturbation $\Delta R(\la,\epsilon)=\epsilon\Delta R(\la) + o(\epsilon)$ satisfying certain technical assumptions (corresponding to perturbing the blocks of a polynomial system matrix), the perturbed rational matrix $R(\la) + \Delta R(\la,\epsilon)$ has a finite zero $\widehat{\la}_0(\epsilon)$ such that $$| \widehat{\la}_0(\epsilon)-\la_0 | = \epsilon \lim_{\la \rightarrow \la_0}\left| \frac{y(\la)^T \Delta R(\la) x(\la)}{y(\la)^T R'(\la) x(\la)} \right| + o(\epsilon).$$

From this result,  after giving some preliminaries on rational matrices and polynomial system matrices in Section \ref{sec:pre}, we define in Section \ref{sec:pert} a structured condition number $\kappa_S$ for a locally minimal \cite{local} polynomial system matrix $P(\la)$ of a rational matrix $R(\la)$, where the structured perturbation $\epsilon \Delta P(\la)$ preserves the degrees of the blocks of $P(\la)$ in \eqref{eq:psm}. For that, we fix the function $f$ to be the input-output function that maps $P(\la)$ to $\la_0$ and $P(\la)+\epsilon \Delta P(\la)$ to $\hat{\la}_0(\epsilon)$. In Section \ref{sec:formula}, we derive an expression for $\kappa_S$. Then, in Section \ref{sec:comparison}, we show that this structured condition number for polynomial system matrices is never larger, and can be much smaller, than Tisseur's unstructured condition number for matrix polynomials \cite{tisseur}. If we consider the case in which the rational matrix $R(\la)$ is expressed in the form $R(\la)=D(\la)+C(\la I_n - A)^{-1}B$, where $D(\la)$ is a polynomial matrix and $C(\la I_n - A)^{-1}B$ is a minimal state-space realization, it was shown in \cite{Backerrors} that there are algorithms to compute the zeros of $R(\la)$, via linearization, that guarantee a small (global) backward error, preserving precisely the structure of $R(\la)$. Note that this is a particular case in the representation \eqref{eq:general representation} where the matrix $A(\la)$ is linear and the matrices $B(\la)$ and $C(\la)$ are constant, but expressing $R(\la)$ in that form is valid for any rational matrix and is used in applications \cite{Rosen70}. Therefore, in problems where the structured condition number is much smaller than the unstructured one, structured algorithms can thus eventually lead to significant increase in accuracy. We also include some numerical experiments on this comparison in Section \ref{sec:num_exp}. Finally, in Section \ref{sec:con} we give some conclusions and related open problems.

\section{Preliminaries}\label{sec:pre}

\subsection{Rational matrices and related subspaces}

A rational matrix $R(\lambda)\in\mathbb{C}(\lambda)^{p\times m}$ is a matrix whose entries are scalar rational functions in the variable $\lambda$. Given $\la_0\in\C$, $R(\lambda)$ is said to be  \textit{defined at} $\la_0$ if $R(\la_0)\in\mathbb{C}^{p\times m}$ and is said to be  \textit{invertible at} $\la_0$ if, in addition, $\det R(\la_0)\neq 0$. Rational matrices can have zeros and poles, which can be defined through the notion of the local Smith--McMillan form. 
That is, for any rational matrix $R(\lambda)\in\mathbb{C}(\lambda)^{p\times m}$, with normal rank $r$, and $\la_0\in\C$ there exist rational matrices $G_{1}(\lambda)$ and $G_{2}(\lambda)$, that are invertible at $\la_0$, such that
	\begin{equation}\label{eq_localSM}
	G_{1}(\lambda)R(\lambda)G_{2}(\lambda)=\left[\begin{array}{cc}
	\text{diag}\left((\la-\la_{0})^{\nu_{1}},\ldots, (\la-\la_{0})^{\nu_{r}}\right)&0 \\
	0& 0_{(p-r)\times (m-r)}
	\end{array}\right],
	\end{equation}
	where $\nu_{1}\leq\cdots\leq\nu_{r}$ are integers. The integers $\nu_{1},\ldots,\nu_{r}$ are uniquely determined by $R(\la)$ and $\la_{0}$ and are called the \textit{invariant orders at} $\la_0$ of $R(\la)$. The matrix in \eqref{eq_localSM} is called the \textit{local Smith--McMillan form of $R(\la)$ at} $\la_{0}.$

 If $\nu_i>0$, for some $i=1,\ldots,r$, then $\la_0$ is a \textit{zero} (or \emph{eigenvalue}) of $R(\la)$ with \textit{partial multiplicity} $\nu_i$, and $(\la-\la_0)^{\nu_i}$ is said to be a \textit{zero elementary divisor of $R(\la)$ at} $\la_0.$ If $\nu_i<0$, for some $i=1,\ldots,r$, then $\la_{0}$ is a \textit{pole} of $R(\la)$ with \textit{partial multiplicity} $-\nu_i$, and $(\la-\la_0)^{-\nu_i}$ is said to be a \textit{pole elementary divisor of $R(\la)$ at} $\la_0$ \cite{Kailath,Vard}. A zero of $R(\la)$ is said to be \textit{simple} if it only has one positive invariant order and it is equal to $1$  (but a simple zero can also be a pole with any pole partial multiplicities). The (global) definition of the Smith--McMillan form of a rational matrix can be found in \cite{Kailath,Rosen70,Vard} and was first given by B. McMillan in \cite{McMi2}.
 
 If $R(\la)$ is a polynomial matrix then the nonzero integers $\nu_i\neq 0$ are all positive and are called \textit{partial multiplicities of $R(\la)$ at} $\la_0$ \cite{GohLR09}. In addition, in this case, the diagonal matrix in \eqref{eq_localSM} is simply called the \textit{local Smith form of $R(\la)$ at $\la_0$} and the polynomials $(\la-\la_0)^{\nu_i}$ with $\nu_i\neq 0$ are called the \textit{elementary divisors of $R(\la)$ at} $\la_0.$ The finite zeros of polynomial matrices are also called eigenvalues. 

For rational matrices having full column normal rank, zeros can be determined by the following equivalent definition \cite[Chapter 27]{notes_Verghese}.

\begin{definition}[Zero location] A rational matrix $R(\lambda)\in\mathbb{C}(\lambda)^{p\times m}$ with full column normal rank has a zero at $\la_0$ if there exists a rational vector $x(\la)\in\mathbb{C}(\lambda)^{m}$ defined at $\la_0$ such that $x(\la_0)\neq 0$ and 
	\begin{equation}\label{eq:lim}
	\lim_{\lambda \to \lambda_0}R(\lambda)x(\la) = 0.
	\end{equation}
\end{definition}

It follows that, in the regular case, i.e., $p=m$ and $\det R(\lambda)\not\equiv 0$, $R(\la)$ has a zero at $\la_0$ if there exist rational vectors $x(\la)\in\mathbb{C}(\lambda)^{m}$ and $y(\la)\in\mathbb{C}(\lambda)^{m}$ defined at $\la_0$ such that 
	\begin{equation}\label{eq:REPlim}
\lim_{\lambda \to \lambda_0}R(\lambda)x(\la) = 0 \quad \quad \mbox{and} \quad \quad \lim_{\lambda \to \lambda_0} y(\la)^T R(\lambda)=0,
\end{equation}
with $x(\la_0)\neq 0$ and $y(\la_0)\neq 0$.
Then $x(\la)$ and $y(\la)$ are said to be right and left \textit{root vectors} of $R(\lambda)$ associated with $\la_0$, respectively \cite{vannof}. Their evaluation at $\la_0$, that is, $x(\la_0)$ and $y(\la_0)$, are called right and left \emph{eigenvectors} of $R(\la)$ associated with $\la_0$.

\begin{remark}\rm The REP is often defined in the literature as the problem of finding scalars $\lambda_0\in\mathbb{C}$ and nonzero vectors $x\in\mathbb{C}^m$ and $y\in\mathbb{C}^m$ such that 
	$R(\lambda_0)x = 0$ and $ y^T R(\lambda_0)=0.$
	But this definition assumes that $\la_0$ is a zero of $R(\la)$ but not a pole, as $R(\la)$ would be defined at $\la_0$. In this paper we consider the more general definition in \eqref{eq:REPlim} to define the REP in order to capture all the zeros of $R(\la)$.
\end{remark}

\begin{remark}\rm
The rational vectors $x(\la)$ and $y(\lambda)$ appearing in \eqref{eq:REPlim} are root vectors of $R(\la)$ for the special case of a regular rational matrix $R(\la)$, see \cite[Definition 3.4]{vannof} for a more general definition of root vectors for rational matrices, not necessarily regular or square, and over arbitray fields. At least in the regular case, root vectors had been studied in the more general context of nonlinear eigenvalue problems over the complex field: see \cite{GT} and the references therein; for the (more involved) nonregular case, a thorough analysis has recently been made for polynomial \cite{DN20} and rational \cite{vannof} matrices over an arbitrary field.

In particular, the theory of root vectors for rational functions can be refined \cite{vannof} so that they not only identify zeros, but also their partial multiplicities (via the concept of a maximal set). As in this paper we will restrict to simple zeros, these details are not important. We refer the reader who would like to know more about maximal sets of root vectors to \cite{DN20,GT,vannof} and the references therein.
\end{remark}

\begin{example}\label{ex_rootvector} Consider the rational matrix
	$R(\la)=\left[\begin{smallmatrix}
	1 & 0 \\
	\frac{1}{\la-1} & 1
	\end{smallmatrix}\right].$ 	It is clear that $R(\la)$ has a pole at $\la_0=1$, as $R(\la)$ is not defined at $1$, but it is less clear that $R(\la)$ has also a zero at $1$. If we consider $x(\la)=\left[\begin{smallmatrix}
\la-1 \\ -1
	\end{smallmatrix}\right]$ we have that $\displaystyle\lim_{\lambda \to 1}R(\lambda)x(\la) = 0$ and $x(1)\neq 0$, i.e., $x(\la)$ is a right root vector of $R(\la)$ associated with $\la_0=1$. Thus $\la_0=1$ is also a zero of $R(\la)$. Indeed, the Smith--McMillan form of $R(\la)$ is $\left[\begin{smallmatrix}
	\frac{1}{\la-1}  & 0 \\
	0 & \la -1
	\end{smallmatrix}\right]$, which shows that $R(\lambda)$ has a zero and a pole at $\lambda_0=1$.

\end{example}

\subsection{Null spaces}

For any rational matrix $R(\la)$ (regular or singular), we define the following rational vector spaces over $\mathbb{C}(\la)$:
\begin{equation*}\label{eq:rat_subspaces}
\begin{array}{l}
\ker R(\la)=\{x(\la)\in\mathbb{C}(\la)^{m\times 1}: R(\la)x(\la)=0 \}, \text{ and}\\
\coker R(\la)=\{y(\la)^T \in\mathbb{C}(\la)^{1\times p}:  y(\la)^T R(\la)=0\},

\end{array}
\end{equation*}
which are called the right and left null spaces over $\mathbb{C}(\la)$ of $R(\la)$, respectively. 

For a finite $\lambda_{0}\in\mathbb{C}$ that is not a pole of $R(\la)$, $\ker R(\la_{0})$ and $\coker R(\la_{0})$ denote the right and left null spaces over $\mathbb{C}$ of the constant matrix $R(\la_{0})$, respectively. 
	Namely,
	
	\[
	\begin{array}{l}
	\ker R(\la_0)=\{x\in\mathbb{C}^{m\times 1}: R(\la_0)x=0\}, \text{ and}\\
	\coker R(\la_0)=\{y^T\in\mathbb{C}^{1\times p}: y^T R(\la_0)=0\}.
	\end{array}
	\]
The vector spaces $\ker R(\la_0)$ and $\coker R(\la_0)$ are called the right and left null spaces of $R(\la)$ at $\la_0$, respectively.

The definition of $\ker R(\la_0)$ and $\coker R(\la_0)$ can be generalized to the case where $\la_0$ is possibly a pole of $R(\la)$ as follows:
	\begin{equation*}
\begin{split}
		\ker R(\la_0):=\{x\in\mathbb{C}^{m\times 1}: & \text{ there exists } x(\la) \in \C(\la)^{m\times 1} \text{ with } \\ &  \displaystyle\lim_{\la \rightarrow \la_0} R(\la)x(\la)=0 \ \mathrm{and} \ x(\la_0)=x\}, \ \mathrm{and}
\end{split}
	\end{equation*}
\begin{equation*}
\begin{split}
	\coker R(\la_0):=\{y^T\in\mathbb{C}^{1\times p}: & \text{ there exists } y(\la)^T \in \C(\la)^{1\times p} \text{ with } \\ &  \displaystyle\lim_{\la \rightarrow \la_0} y(\la)^T R(\la)=0 \ \mathrm{and} \ y(\la_0)=y\}.
	\end{split}
	\end{equation*}
	
	For regular rational matrices $R(\la)$,	if $\la_0$ is a zero of $R(\la)$ then $\ker R(\la_{0})$ and $\coker R(\la_{0})$ are non trivial and contain the right and left eigenvectors of $R(\la)$ associated with $\la_0$, respectively. In \cite[Definition 3.8 and Remark 3.9]{vannof} one can find a more general definition for $\ker R(\la_{0})$ and $\coker R(\la_{0})$ valid over an arbitrary field $\mathbb{F}$ and that specializes to the one above when $\mathbb{F}=\mathbb{C}$. 
	

\subsection{Polynomial system matrices and transfer function matrices}

It is known that any rational matrix $R(\lambda)\in\mathbb{C}(\lambda)^{p\times m}$ can be written with an expression of the form
\begin{equation}\label{eq:rational matrix}
R(\lambda)= D(\lambda) + C(\lambda)A(\lambda)^{-1}B(\lambda),
\end{equation}
where $A(\la)$, $B(\la)$, $C(\la)$ and $D(\la)$ are  matrix polynomials, with $A(\lambda)\in\mathbb{C}[\la]^ {n\times n}$ regular. In addition, the representation \eqref{eq:rational matrix} is not unique. The associated polynomial matrix
\begin{equation}\label{eq_polsysmat}
P(\la)=\begin{bmatrix}
-A(\la) & B(\la)\\
C(\la) & D(\la)
\end{bmatrix}
\end{equation}
is said to be a polynomial system matrix of the rational matrix $R(\la)$ in \eqref{eq:rational matrix}, and the rational matrix $R(\la)$ is called the transfer function matrix of $P(\la)$ (see \cite{Rosen70}). Note that $R(\la)$ is the Schur complement of $A(\la)$ in $P(\la)$. The polynomial system matrix $P(\la)$ in \eqref{eq_polsysmat}, with $n>0$, is said to be minimal if the polynomial matrices 
\begin{equation} \label{mat_minimalidad}
\begin{bmatrix} -A(\la) \\ C(\la)\end{bmatrix}\quad \text{and} \quad 	\begin{bmatrix} A(\la) & B(\la) \end{bmatrix}
\end{equation}
have no eigenvalues in $\mathbb{C}$ \cite{Rosen70}. An important consequence of minimality is that, if a polynomial system matrix $P(\la)$ is minimal, the poles of its transfer function matrix $R(\la)$ are the zeros of $A(\la)$, and the zeros of $R(\la)$ are the zeros of $P(\la)$, together with their partial multiplicities \cite[Chapter 3, Theorem 4.1]{Rosen70}. The notion of minimality was extended in a local sense in \cite{local}. For that, if the matrices in \eqref{mat_minimalidad} have no eigenvalues in a nonempty subset $\Omega\subseteq\mathbb{C}$, $P(\la)$ is said to be minimal in $\Omega$ \cite[Definition 3.1]{local}. Then, if a polynomial system matrix is minimal in a set $\Omega$, the poles of $R(\la)$ in $\Omega$ are the zeros of $A(\la)$ in $\Omega$, and the zeros of $R(\la)$ in $\Omega$ are the zeros of $P(\la)$ in $\Omega$, together with their partial multiplicities \cite[Theorem 3.5]{local}. If $P(\la)$ is minimal in $\Omega:=\{\la_0\}$, with $\la_0\in\C$, then $P(\la)$ is said to be minimal at $\la_0$.

Given a polynomial system matrix $P(\la)$ minimal at $\la_0$, the following Proposition \ref{prop:eigenvectorslimit} establishes the relation between the right and left null spaces of $P(\la)$ and its transfer function matrix $R(\la)$ at $\la_0$, and it is valid for arbitrary rational matrices $R(\la)$ (regular or singular). This result was stated in \cite[Lemma 2.1]{BFR} in the (easier) special case when $\la_0$ is not a pole of $R(\la)$, which follows from \cite[Proposition 5.1]{DMQ} and \cite[Proposition 5.2]{DMQ}. Here, we give a proof of the general case.

\begin{proposition}\label{prop:eigenvectorslimit}
	
		Let 
	\[
	P(\la)=\begin{bmatrix}
	-A(\la) & B(\la)\\
	C(\la) & D(\la)
	\end{bmatrix}\in \mathbb{C}[\la]^{(n+p)\times (n+m)}
	\]
	be a polynomial system matrix, with $A(\la)\in\mathbb{C}[\lambda]^{n\times n}$ regular, and transfer function matrix $R(\la) = D(\lambda)+C(\lambda)A(\lambda)^{-1}B(\lambda)\in \mathbb{C}(\la)^{p\times m}$. Let $\la_0\in\mathbb{C}$. If $P(\la)$ is minimal at $\la_0$ then the following statements hold:
	\begin{itemize}
		\item[\rm(a)] The linear map
		\begin{equation*}
		L_r \, : \, \ker R(\la_0) \longrightarrow \ker P(\la_0), \quad
		x_0 \longmapsto   \begin{bmatrix} \displaystyle\lim_{\lambda \to \lambda_0} A(\lambda)^{-1}B(\lambda)x(\la)\\x_0 \end{bmatrix}
		\end{equation*}
		is a bijection, where $x(\la)\in\C(\la)^m$ is such that $\displaystyle\lim_{\lambda \to \lambda_0}R(\lambda)x(\la) = 0$ and $x(\la_0)=x_0.$
		\item[{\rm(b)}] The linear map
		\begin{equation*}
		L_\ell \, : \, \coker R(\la_0)  \longrightarrow \coker P(\la_0) , \quad
		y_0^T  \longmapsto  \begin{bmatrix}  \displaystyle\lim_{\lambda \to \lambda_0} y(\la)^TC(\lambda)A(\lambda)^{-1} & y_0^T \end{bmatrix}
		\end{equation*}
		is a bijection, where $y(\la)\in\C(\la)^p$ is such that $\displaystyle\lim_{\lambda \to \lambda_0}y(\la)^TR(\lambda) = 0$ and $y(\la_0)=y_0.$
	\end{itemize}

\end{proposition}

\begin{proof} We only prove part $\rm (a)$ since part $\rm (b) $ is analogous.

 First, we prove that $L_r$ is well-defined.
 Let $x_0\in \ker R(\la_0) $. Then there exists a rational vector $x(\la)\in\C(\la)^m$ such that $\displaystyle\lim_{\lambda \to \lambda_0}R(\lambda)x(\la) = 0$ and $x(\la_0)=x_0$, by definition. Since $\displaystyle\lim_{\lambda \to \lambda_0}R(\lambda)x(\la) = 0$ we have that
	\begin{equation}\label{eq:limpolsystem}
\displaystyle\lim_{\lambda \to \lambda_0}\begin{bmatrix}
-A(\la) & B(\la)\\
C(\la) & D(\la)
\end{bmatrix}\begin{bmatrix} A(\lambda)^{-1}B(\lambda)x(\la)\\x(\la) \end{bmatrix}=0.
	\end{equation}
Then, note that \eqref{eq:limpolsystem} can be rewritten as
	\begin{equation}\label{eq_limitminimal}
\displaystyle\lim_{\lambda \to \lambda_0}\left(\begin{bmatrix}
-A(\la) \\
C(\la) 
\end{bmatrix}A(\lambda)^{-1}B(\lambda)x(\la) + \begin{bmatrix} B(\lambda) \\ D(\la)\end{bmatrix}x(\la)\right)=0.
	\end{equation}
Since $P(\la)$ is minimal at $\la_0$, we know that $\left[\begin{smallmatrix}
-A(\la) \\
C(\la) 
\end{smallmatrix}\right]$ has a rational left inverse that is defined at $\la_0$. 
In other words, there exists a rational matrix $H(\la)$ such that $\la_0$ is not a pole of $H(\la)$ and $H(\la)\left[\begin{smallmatrix}
-A(\la) \\
C(\la) 
\end{smallmatrix}\right]=I$ (the claim follows immediately from the observation that the local Smith form of $\left[\begin{smallmatrix}
-A(\la) \\
C(\la) 
\end{smallmatrix}\right]$ at $\la_0$ is $\left[\begin{smallmatrix}
I \\
0
\end{smallmatrix}\right]$). Taking into account \eqref{eq_limitminimal}, we obtain that the limit
\begin{align*}\label{eq_limitminimal2}
\displaystyle\lim_{\lambda \to \lambda_0}A(\lambda)^{-1}B(\lambda)x(\la)
&=
\displaystyle\lim_{\lambda \to \lambda_0}H(\la)\begin{bmatrix}
-A(\la) \\
C(\la) 
\end{bmatrix}A(\lambda)^{-1}B(\lambda)x(\la)
  \\
&= -\lim_{\la\to\la_0} H(\la) \begin{bmatrix} B(\la) \\ D(\la) \end{bmatrix} x(\lambda) =
H(\la_0) \begin{bmatrix} B(\la_0) \\ D(\la_0) \end{bmatrix} x
\end{align*}
is defined. 
Thus, $\displaystyle\lim_{\lambda \to \lambda_0}\left[\begin{smallmatrix} A(\lambda)^{-1}B(\lambda)x(\la)\\x(\la) \end{smallmatrix}\right]\in \ker P(\la_0)$. This proves that $L_r$ is well-defined.  
Next, note that if $\left[\begin{smallmatrix} \lim_{\lambda \to \lambda_0} A(\lambda)^{-1}B(\lambda)x(\la)\\x_0 \end{smallmatrix}\right]=0$ then $x_0=0$, and this establishes the injectivity of $L_r$. 
Finally, to see that $L_r$ is a bijection, we prove that $\dim \ker R(\la_0) = \dim \ker P(\la_0) $. Denote by $s$ the number of positive invariant orders of $R(\la)$ at $\la_0$; observe that $s$ is also equal to the number of partial multiplicities of $P(\la)$ at $\la_0$, by \cite[Theorem 3.5]{local}, since $P(\la)$ is minimal at $\la_0$.  By \cite[Theorem 3.10]{vannof},
$$\dim \ker R(\la_0) = \dim \ker R(\la)+s, \qquad \dim \ker P(\la_0) = \dim \ker P(\la)+s.$$
On the other hand, by \cite[Lemma 2.1]{PQ}, $\dim \ker P(\la) = \dim \ker R(\la)$, and hence $\dim \ker P(\la_0) = \dim \ker R(\la_0)$. 
\end{proof}

\begin{remark}[Notation] \rm \label{rem:notation} If $\la_0$ is a pole of $R(\la)$, then $A(\la)^{-1}$ is not defined at $\la_0$. However, as we proved in Proposition \ref{prop:eigenvectorslimit}, if $P(\la)$ is minimal at $\la_0$ the limit $\displaystyle\lim_{\lambda \to \lambda_0}A(\lambda)^{-1}B(\lambda)x(\la)$ is defined at $\la_0$, even if $\la_0$ is a pole of $R(\la)$. Thus, in what follows, for a given $x(\la)\in \mathbb{C}(\la)^m$ such that $\displaystyle\lim_{\lambda \to \lambda_0}R(\lambda)x(\la) = 0$, we will slightly abuse notation and write $A(\lambda_0)^{-1}B(\lambda_0)x(\la_0)$ as a shorthand for $\displaystyle\lim_{\lambda \to \lambda_0}A(\lambda)^{-1}B(\lambda)x(\la)$; this should be interpreted as the rational vector $A(\lambda)^{-1}B(\lambda)x(\la)$ evaluated at $\la_0$ (this evualuation is guaranteed to yield a complex vector by Proposition \ref{prop:eigenvectorslimit}). The same notation will be used for the limit $\displaystyle\lim_{\lambda \to \lambda_0}y(\la)^TC(\la)A(\lambda)^{-1}$ for any $y(\la)\in\C(\la)^p$ such that $\displaystyle\lim_{\lambda \to \lambda_0}y(\la)^TR(\lambda) = 0$, i.e., $y(\la_0)^TC(\la_0)A(\lambda_0)^{-1}$.
	
\end{remark}

We know that if $\la_0$ is a zero of a rational matrix $R(\la)$ then $\la_0$ is also a zero (an eigenvalue) of any polynomial system matrix $P(\la)$ of $R(\la)$, with same partial multiplicities, whenever $P(\la)$ is minimal at $\la_0$ \cite[Theorem 3.5]{local}. Then Proposition \ref{prop:eigenvectorslimit} shows a corresponding relation between the eigenvectors, even if $\la_0$ is not only a zero but also a pole of $R(\la)$. Moreover, since the maps in Proposition \ref{prop:eigenvectorslimit} are bijections, they preserve linear independence. Thus, one can recover a basis of $\ker P(\la_0)$ (resp. $\coker P(\la_0)$) from a basis of $\ker R(\la_0)$ (resp. $\coker R(\la_0)$), and conversely.

\begin{example} Consider the rational matrix $R(\la)$ in Example \ref{ex_rootvector}, and the right root vector $x(\la)$ of $R(\la)$ associated with $\la_0=1$. Define the polynomial system matrix
	$$P(\la):=\left[\begin{array}{c|cc}
\la-1 & 1 & 0 \\ \hline
0 & 1 & 0 \\
-1 & 0 & 1
	\end{array}\right] =:\begin{bmatrix}
	-A(\la) & B(\la)\\
	C(\la) & D(\la)
	\end{bmatrix},$$
	whose transfer function matrix is $R(\la)$ and is minimal at $\la_0=1$. Then $\la_0=1$ is also a zero of $P(\la)$ and $v=\left[\begin{smallmatrix}  A(\lambda_0)^{-1}B(\lambda_0)x(\la_0)\\x(\la_0) \end{smallmatrix}\right]=\begin{bmatrix} -1 & 0 & -1 \end{bmatrix}^T$ is a right eigenvector of $P(\la)$ associated with $\la_0=1$.

\end{example}

\section{Perturbation expansion for simple zeros}\label{sec:pert}
Let $R(\la)$ be the transfer function matrix of a polynomial system matrix $P(\la)$. In this section, we assume that $R(\la)$ is regular (so that $P(\la)$ is also regular) and that $\lambda_0$ is a simple finite zero of $R(\la)$. If $P(\la)$ is minimal at $\la_0$, $\la_0$ is also a simple eigenvalue of $P(\la)$ \cite[Theorem 3.5]{local}. We first study how analytic perturbations of $P(\la)$ affect to simple eigenvalues and local minimality.  
\begin{theorem}\label{th:analytic}
Let $ P(\lambda) =\left[\begin{smallmatrix}
-A(\lambda) & B(\lambda)\\
C(\lambda) & D(\lambda)
\end{smallmatrix}\right]$ be a regular polynomial system matrix, with $A(\la)$ regular. Let $\lambda_0\in\C$ be a finite simple eigenvalue of $P(\la)$, with right (resp. left) eigenvector $v$ (resp. $w$), and assume that $P(\la)$ is minimal at $\la_0$. Let
 $$
 \widehat P(\lambda,\epsilon)=\begin{bmatrix}
-\widehat A(\lambda,\epsilon) & \widehat B(\lambda,\epsilon)\\
\widehat C(\lambda,\epsilon) & \widehat D(\lambda,\epsilon)
\end{bmatrix}
$$
be a matrix-valued function that can be written as $\widehat P(\lambda,\epsilon)=\sum_{i=0}^d P_i(\epsilon) \lambda^i$ where each $P_i(\epsilon)$ is an analytic matrix function of $\epsilon$ in some open ball centered at $0$. Assume that $\widehat P(\lambda,0)=P(\lambda)$
Then there exists a constant $r>0$ such that the following statements hold: 
	\begin{itemize}
		\item[\rm (a)] There exist functions $\widehat \lambda_0(\epsilon)$, $\widehat{v}(\epsilon)$ and $\widehat{w}(\epsilon)$ such that {\rm (i)} $\widehat \lambda_0(\epsilon), \widehat{v}(\epsilon)$ and $ \widehat{w}(\epsilon)$ are analytic in $\epsilon$ for all $|\epsilon|<r$; and {\rm (ii)} $\widehat \lambda_0(0)=\lambda_0$, $\widehat{w}(0)=w$ and $ \widehat{v}(0)=v$.
		\item[\rm (b)] 	
		$ \widehat P(\la,\epsilon)$ has a simple zero at $\widehat \lambda_0(\epsilon)$, for all $|\epsilon|<r$, with associated right (resp. left) eigenvector $\widehat{v}(\epsilon)$ (resp. $\widehat{w}(\epsilon)$).
		\item[\rm (c)] $\widehat P(\lambda,\epsilon)$, viewed as a polynomial matrix in $\la$, is a polynomial system matrix minimal at both $\la_0$ and $\widehat \lambda_0(\epsilon)$ for all $|\epsilon|<r$.
		\item[\rm (d)] The transfer function matrix $ \widehat R(\la,\epsilon)$ of $ \widehat P(\la,\epsilon)$ has a simple zero at $\widehat \lambda_0(\epsilon)$ for all $|\epsilon|<r$.
	\end{itemize}
\end{theorem}
\begin{proof}
Items (a) and (b) are classical results, see for example \cite{ACL93,Kato}.

For $\rm(c)$, we only prove that for sufficiently small $\epsilon$, $\widehat P(\lambda,\epsilon)$ must be minimal at $\widehat{\lambda}_0(\epsilon)$. Indeed, minimality is equivalent to $\begin{bmatrix}
-\widehat{A}(\widehat{\lambda}_0(\epsilon),\epsilon)&\widehat{B}(\widehat{\lambda}_0(\epsilon),\epsilon)
\end{bmatrix}$ having full row rank and $\begin{bmatrix}
-\widehat{A}(\widehat{\lambda}_0(\epsilon),\epsilon)\\\widehat{C}(\widehat{\lambda}_0(\epsilon),\epsilon)
\end{bmatrix}$ having full column rank. On the other hand, $\widehat P(\widehat{\lambda}_0(\epsilon),\epsilon)$ is analytic in $\epsilon$; hence, so is any of its submatrices and the statement follows by the continuity of singular values \cite{Kato}.

Finally, $\rm(d)$ follows from $\rm (b)$ and $\rm (c)$. That is, since $\widehat P(\la,\epsilon)$ is minimal at $\widehat{\lambda}_0(\epsilon)$ and has a simple zero at $\widehat \lambda_0(\epsilon)$, by \cite[Theorem 3.5]{local}, it implies that $\widehat{\la}_0(\epsilon)$ is also a simple zero of $\widehat{R}(\lambda,\epsilon)$.
\end{proof}

 \begin{remark} \rm For a given rational matrix $R(\la)$, not every small perturbation of $R(\la)$ has a polynomial system matrix $\widehat P(\lambda,\epsilon)$ satisfying the assumptions in Theorem \ref{th:analytic}. For instance, consider the rational matrix
		$$ R(\lambda)= \begin{bmatrix}
		\frac{1}{\lambda} & \frac{1}{\lambda^3}\\
		0 &\frac{1}{\lambda}
		\end{bmatrix}.$$ 
		The invariant orders of $R(\la)$ at $0$ are $-3,1$, and there are not any other finite zeros or poles. For any $\epsilon >0$, let
		$$ R(\lambda)+\epsilon \Delta R(\lambda) = \frac{1}{\lambda^3}  \begin{bmatrix}
		\lambda^2 + \epsilon & 1\\
		0 & \lambda^2-\epsilon
		\end{bmatrix},\quad \text{where} \quad \Delta R(\lambda)=\begin{bmatrix}
		\frac{1}{\lambda^3} & 0\\
		0 & -\frac{1}{\lambda^3}
		\end{bmatrix}.$$
		The invariant orders of $R(\lambda)+\epsilon \Delta R(\lambda)$ at $0$ are $-3,-3$. Then, $0$ is a pole. Moreover, the zeros of $R(\lambda)+\epsilon\Delta R(\lambda)$ are \emph{not} analytic in $\epsilon$: indeed, they can be shown to be equal to $\pm \epsilon^{1/2}$ and $\pm i \epsilon^{1/2}$. Thus, the statement of Theorem \ref{th:analytic}$\rm (d)$ is not satisfied for this particular perturbation. However, this is not a counterexample to the theorem because the assumptions are also false: suppose indeed that $P(\lambda)=\widehat P(\lambda,0)$ such that $P(\la)$ and $\widehat P(\la,\epsilon)$ are polynomial system matrices of $R(\la)$ and $R(\la)+\epsilon \Delta R(\la)$, respectively, both minimal at $0$. Then, by \cite[Theorem 3.5]{local}, the matrices $A(\la)$ and $\widehat{A}(\la,\epsilon)$ must have determinant, respectively, $p(\la) \lambda^3$ and $f(\lambda,\epsilon) \lambda^6$ where $p(\la)$ is a polynomial with $p(0)\neq 0$ and $f(\lambda,\epsilon)$ is a function polynomial in $\la$, analytic in $\epsilon$ with $f(0,\epsilon)\neq 0.$ But then, if $P(\lambda)=\widehat P(\lambda,0)$, by the continuity of the determinant we would have that $f(\lambda,0)\lambda^6=p(\la) \lambda^3$, which is impossible. Although every rational matrix admits minimal polynomial system matrices, the crux here is the requirement $\widehat P(\la,0)=P(\la)$. For example, $$P(\la)=\left[\begin{array}{c|cc}
		\la^3 & \la^2 & 1 \\\hline
		1 & 0 & 0 \\
		\la^2 & -\la & 0 
		\end{array}\right] \quad\text{and} \quad Q(\la,\epsilon)=\left[\begin{array}{cc|cc}
		\la^3 & 0 & \epsilon & 0\\
		0 & \la^3 & \la^2 & 1 \\ \hline
		0 & 1 & 0 & 0 \\
		\epsilon &	\la^2 & -\la & 0 
		\end{array}\right]$$
		are minimal polynomial system matrices of $R(\lambda)$ and $R(\la)+\epsilon \Delta R(\la)$, respectively, but clearly $Q(\lambda,0)\neq P(\la)$ (as the sizes are different).
		
	\end{remark}

Nevertheless, in this paper we consider perturbed rational matrices having polynomial system matrices that always satisfy the assumptions in Theorem \ref{th:analytic}. More precisely, we consider rational matrices $R(\la)$ expressed as in \eqref{eq:rational matrix}, i.e., $R(\lambda)= D(\lambda) + C(\lambda)A(\lambda)^{-1}B(\lambda)$, for arbitrary polynomial matrices $A(\la)$, $B(\la)$, $C(\la)$ and $D(\la)$, and its associated polynomial system matrix $P(\la)=\left[\begin{smallmatrix}
-A(\la) & B(\la)\\
C(\la) & D(\la)
\end{smallmatrix}\right].$ Then, we consider a perturbation of $P(\lambda)$ of the form:
\begin{equation}\label{eq_polsysmat_per}
P(\la)+\epsilon \Delta P(\lambda):=\begin{bmatrix}
-A(\la) & B(\la)\\
C(\la) & D(\la)
\end{bmatrix}+ \epsilon \begin{bmatrix}
-\Delta A(\la) & \Delta B(\la)\\
 \Delta C(\la) & \Delta D(\la)
\end{bmatrix},
\end{equation} 
 where $\Delta A(\la)$, $\Delta B(\la)$, $\Delta C(\la)$ and $\Delta D(\la)$ are polynomial matrices. It is clear that $ \widehat P(\lambda,\epsilon):=P(\la)+\epsilon \Delta P(\lambda)$ satisfies the assumptions in Theorem \ref{th:analytic} as it is a linear perturbation in $\epsilon$ of $P(\la)$. In order to study perturbations of simple eigenvalues, we will use Lemma \ref{lem_der}, which is known for general polynomial matrices \cite[Proof of Theorem 5]{tisseur} but that we state here for the particular case of polynomial system matrices.

\begin{lemma}\label{lem_der} Let $\la_0\in\C$ be a simple eigenvalue of a regular polynomial system matrix $P(\la)$; and let $v$ and $w$ be, respectively, right and left eigenvectors of $P(\lambda)$ associated with $\lambda_0$. Consider the perturbed polynomial system matrix matrix $P(\la)+\epsilon \Delta P(\lambda)$.
	Then, for $\epsilon>0$ small enough, there exist $\Delta \lambda_0$ and $\Delta v$ satisfying 
$$\left[P(\la_0+ \Delta\la_0)+\epsilon\Delta P(\la_0+\Delta \la_0)\right][v+\epsilon\Delta v]=0.$$  That is, $\la_0 +  \Delta \la_0$ is an eigenvalue for the perturbed polynomial system matrix $P(\la)+\epsilon\Delta P(\la)$. In addition, 
\begin{equation}\label{eq_pertzero22}
|\Delta \lambda_0| \,=\, \epsilon \left| \dfrac{w^T\Delta P(\la_0) v}{w^T P^\prime (\la_0) v } \right|+o(\epsilon).
\end{equation}
\end{lemma}


Now, we consider the perturbed transfer function matrix of the perturbed polynomial system matrix in \eqref{eq_polsysmat_per}. Namely,
\begin{align}\label{eq:ptfm}
\begin{split}
R(\lambda) + \Delta R(\lambda,\epsilon):= D(\lambda)+\epsilon & \Delta D(\lambda)+\\
&(C(\lambda)+\epsilon\Delta C(\lambda))(A(\lambda)+\epsilon\Delta A(\lambda))^{-1}(B(\lambda)+\epsilon\Delta B(\lambda)).
\end{split}
\end{align}


Then, we set  
\begin{equation}\label{eq_delta2}
\begin{split}
\Delta R(\lambda) := \left[ \frac{\partial \Delta R(\la,\epsilon)}{\partial \epsilon}\right]_{\epsilon=0} =  \Delta D (\la)&+\Delta C(\la) A(\lambda)^{-1}B(\lambda)+   C(\lambda)A(\lambda)^{-1} \Delta B(\la) \\ 
& - C(\lambda)A(\lambda)^{-1} \Delta A(\la)   A(\lambda)^{-1} B(\lambda).
\end{split}
\end{equation}
Note that \eqref{eq_delta2} implies the formal \footnote{While one can always write such a formal expansion, when $\la \rightarrow \la_0$ and $\la_0$ is a pole of $R(\la)$ then the radius of convergence (in $\epsilon$) of such an expansion may tend to $0$. For example, consider $$R(\la)+\Delta R(\la,\epsilon)=\frac{1}{\la-\epsilon} = \frac{1}{\la} \sum_{k=0}^\infty \frac{\epsilon^k}{\la^k}, \text{ where } R(\la)=\frac{1}{\la} \text{ and } \Delta R(\la,\epsilon)=\frac{\epsilon}{\la^2} + o(\epsilon).$$ For any fixed $\la \neq 0$, the radius of convergence of the power series in $\epsilon$ is $|\lambda|$. Hence, there is no uniform radius of convergence for $\la$ in a neighbourhood of the pole $\la_0=0$.} expansion $\Delta R(\la,\epsilon)=\epsilon \Delta R(\la)+o(\epsilon)$.

In Corollary \ref{cor:der_rat}, we write \eqref{eq_pertzero22} in terms of $R(\la)$ and $\Delta R(\lambda)$ by using root vectors \cite{vannof}. For that, we need Lemma \ref{lem:deno_nume}. 

\begin{lemma}\label{lem:deno_nume} 
	Let $P(\la) \in \C[\la]^{(n+p)\times(n+p)}$ 
	be a polynomial system matrix as in \eqref{eq_polsysmat}, with $A(\la)\in\mathbb{C}[\lambda]^{n\times n}$ regular, and transfer function matrix $R(\la) = D(\lambda)+C(\lambda)A(\lambda)^{-1}B(\lambda)\in \mathbb{C}(\la)^{p\times p}$. 
Assume that $P(\la)$ is minimal at $\la_0$. 
Let $v=\left[\begin{smallmatrix}
		v_1\\ v_2
		\end{smallmatrix}\right]\in\ker P(\la_0)$ and $w^T=\begin{bmatrix}
		w_1^T & w_2^T
		\end{bmatrix}\in\coker P(\la_0)$, partitioned conformably to the blocks of $P(\la)$, with $v\neq 0$ and $w\neq 0$. 
Then, there exist right and left root vectors $x(\la)\in\mathbb{C}(\la)^{p}$ and $y(\la)\in\mathbb{C}(\la)^{p}$ of $R(\lambda)$, both associated with $\la_0$,  such that 
\begin{equation}\label{eq_der}
w^T P'(\la_0)v=\displaystyle\lim_{\lambda \to \lambda_0} y(\la)^T R'(\la) x(\la),
\end{equation}
with $x(\la_0) = v_2$ and $y(\la_0) = w_2$, and that 
\begin{equation}
w^T \Delta P(\la_0)v=\displaystyle\lim_{\lambda \to \lambda_0} y(\la)^T \Delta R(\la) x(\la),
\end{equation}
where $\Delta R(\la)$ is defined as in \eqref{eq_delta2}.
\end{lemma}

\begin{proof} 
By Proposition \ref{prop:eigenvectorslimit}, we know that if $v\in\ker P(\la_0 )$ and $w\in\coker P(\la_0)$ then $v=  \begin{bmatrix} \displaystyle\lim_{\lambda \to \lambda_0} A(\lambda)^{-1}B(\lambda)x(\la)\\x(\la_0) \end{bmatrix}$ and $  w^T=  \begin{bmatrix}\displaystyle\lim_{\lambda \to \lambda_0} y(\la)^TC(\lambda)A(\lambda)^{-1}  & y(\la_0)^T \end{bmatrix}$ for some $x(\la)\in\mathbb{C}(\la)^p$ and $y(\la)\in\mathbb{C}(\la)^p$ with $\displaystyle\lim_{\lambda \to \lambda_0}R(\lambda)x(\la) = 0$ and $\displaystyle\lim_{\lambda \to \lambda_0}y(\la)^T R(\lambda) = 0$, respectively. Then	
	{\small	\begin{equation*}
		\begin{split}
		&	w^TP'(\la_0)v \\ & = \begin{bmatrix}\displaystyle\lim_{\lambda \to \lambda_0}y(\la)^TC(\lambda)A(\lambda)^{-1} & y(\la_0)^T \end{bmatrix} \begin{bmatrix}
		-A'(\la_0) & B'(\la_0)\\
		C'(\la_0) & D'(\la_0)
		\end{bmatrix} \begin{bmatrix} \displaystyle\lim_{\lambda \to \lambda_0} A(\lambda)^{-1}B(\lambda)x(\la)\\x(\la_0) \end{bmatrix}\\
		& = \displaystyle\lim_{\lambda \to \lambda_0} y(\la)^T(D'(\la)+C'(\la)A(\la)^{-1}B(\la)+ C(\la)A(\la)^{-1}B'(\la)\\ & \hspace{7cm} -C(\la)A(\la)^{-1}A'(\la)A(\la)^{-1}B(\la))x(\la) \\
&= \displaystyle\lim_{\lambda \to \lambda_0} y(\la)^T R'(\la)x(\la)
		\end{split}
		\end{equation*}}%
and 
{\small\begin{equation*}
	\begin{split}
	&	w^T\Delta P(\la_0)v \\ & = \begin{bmatrix}\displaystyle\lim_{\lambda \to \lambda_0}y(\la)^TC(\lambda)A(\lambda)^{-1} & y(\la_0)^T \end{bmatrix} \begin{bmatrix}
	-\Delta A(\la_0) & \Delta B(\la_0)\\
	\Delta C(\la_0) & \Delta D(\la_0)
	\end{bmatrix} \begin{bmatrix} \displaystyle\lim_{\lambda \to \lambda_0} A(\lambda)^{-1}B(\lambda)x(\la)\\x(\la_0) \end{bmatrix}\\
	& = \displaystyle\lim_{\lambda \to \lambda_0} y(\la)^T(\Delta D(\la)+\Delta C(\la)A(\la)^{-1}B(\la)+ C(\la)A(\la)^{-1}\Delta B(\la)\\ & \hspace{6cm} -C(\la)A(\la)^{-1}\Delta A(\la)A(\la)^{-1}B(\la))x(\la)\\ &
	= \displaystyle\lim_{\lambda \to \lambda_0} y(\la)^T \Delta R(\la)x(\la),
	\end{split}
	\end{equation*}	}%
	with $x(\la_0) = v_2$ and $y(\la_0) = w_2$. We still have to prove that $x(\la_0)\neq 0$ and $y(\la_0)\neq 0$. We only prove that $x(\la_0)\neq 0$, since $y(\la_0)\neq 0$ can be proved analogously. Since $\displaystyle\lim_{\lambda \to \lambda_0}R(\lambda)x(\la) = 0$ we have that
		\begin{equation*}
		\displaystyle\lim_{\lambda \to \lambda_0}\begin{bmatrix}
		-A(\la) \\
		C(\la) 
		\end{bmatrix}A(\lambda)^{-1}B(\lambda)x(\la) + \begin{bmatrix} B(\lambda) \\ D(\la)\end{bmatrix}x(\la)=0.
		\end{equation*}
		By contradiction, if we assume that $x(\la_0)=0$, we get that $\displaystyle\lim_{\lambda \to \lambda_0} \left[\begin{smallmatrix} B(\lambda) \\ D(\la)\end{smallmatrix}\right] x(\la)=0$ and, therefore, 
		$
		\displaystyle\lim_{\lambda \to \lambda_0}\left[\begin{smallmatrix}
		-A(\la) \\
		C(\la) 
		\end{smallmatrix}\right]A(\lambda)^{-1}B(\lambda)x(\la) =0.
		$
		Taking into account that $P(\la)$ is minimal at $\la_0$, this is, that $\displaystyle\lim_{\la\to\la_0}\left[\begin{smallmatrix}-A(\la) \\ C(\la)\end{smallmatrix}\right] = 
\left[\begin{smallmatrix}-A(\la_0) \\ C(\la_0)\end{smallmatrix}\right]$ has full column rank, we obtain that
		$\displaystyle\lim_{\lambda \to \lambda_0} A(\lambda)^{-1}B(\lambda)x(\la)= 0$, which is a contradiction since this implies that $v=0$. 
\end{proof}

\begin{remark} \rm We remark that $w^T P'(\la_0)v\neq 0$ by \cite[Theorem 3.2]{ACL93}, for right and left eigenvectors $v$ and $w$ of $P(\la)$, respectively, associated with a simple eigenvalue $\la_0$. Taking into account Lemma \ref{lem:deno_nume}, this also implies that $\displaystyle\lim_{\lambda \to \lambda_0} y(\la)^T R'(\la) x(\la)\neq 0,$ for right and left root vectors $x(\la)$ and $y(\la)$ of $R(\la)$, respectively, associated with a simple zero $\la_0$ of $R(\la)$.
\end{remark}

\begin{corollary}\label{cor:der_rat} Assume that the polynomial system matrix $P(\la)$ in Lemma \ref{lem_der} is minimal at $\la_0$, and let $R(\la)$ be the transfer function matrix of $P(\la)$. Let $R(\la)+\Delta R(\lambda,\epsilon)$ be the transfer function matrix of the perturbed polynomial system matrix $P(\la)+\epsilon \Delta P(\lambda)$ in \eqref{eq_polsysmat_per}. Define $\Delta R(\la)$ as in \eqref{eq_delta2}. Then, for $\epsilon>0$ small enough, $\la_0  + \Delta \la_0 $ is also a zero for $R(\la)+\Delta R(\lambda,\epsilon)$. In addition, there exist right and left root vectors $x(\la)\in\mathbb{C}(\la)^{p}$ and $y(\la)\in\mathbb{C}(\la)^{p}$ of $R(\lambda)$ associated with $\la_0$, respectively, such that 
	\begin{equation}\label{eq_pertzero}
	|\Delta \lambda_0| \,=\, \epsilon \lim_{\la\rightarrow \lambda_0}  \left| \dfrac{y(\la)^T\Delta R(\la) x(\la)}{y(\la)^T R^\prime (\la) x(\la) } \right|+o(\epsilon).
	\end{equation}

\end{corollary}

\section{Condition number formulae}\label{sec:formula}

In this section, we define absolute condition numbers for simple zeros of minimal polynomial system matrices as in \eqref{eq:psm}. For that, we express the block
polynomial matrices $A(\lambda)$, $B(\lambda)$, $C(\lambda)$ and $D(\lambda)$ in terms of the monomial basis. Namely,
 \begin{align*}
 &D(\lambda)=\sum_{i=0}^{k_D} D_i\lambda^i\in\mathbb{C}[\lambda]^{p\times p},\quad C(\lambda)=\sum_{i=0}^{k_C} C_i\lambda^i\in\mathbb{C}[\lambda]^{p\times n},\\
 &A(\lambda)=\sum_{i=0}^{k_A} A_i\lambda^i\in\mathbb{C}[\lambda]^{n\times n}, \quad \mbox{and}\quad
 B(\lambda)=\sum_{i=0}^{k_B} B_i\lambda^i\in\mathbb{C}[\lambda]^{n\times p}.
 \end{align*} 
 Then, we perturb the matrix coefficients of each matrix polynomial, and the perturbed rational matrix we get is \eqref{eq:ptfm} with 
\begin{align*}
&\Delta D(\lambda)=\sum_{i=0}^{k_D} \Delta D_i\lambda^i\in\mathbb{C}[\lambda]^{p\times p},\quad \Delta C(\lambda)=\sum_{i=0}^{k_C} \Delta C_i\lambda^i\in\mathbb{C}[\lambda]^{p\times n},\\
&\Delta A(\lambda)=\sum_{i=0}^{k_A} \Delta A_i\lambda^i\in\mathbb{C}[\lambda]^{n\times n}, \quad \mbox{and}\quad \Delta B(\lambda)=\sum_{i=0}^{k_B} \Delta B_i\lambda^i\in\mathbb{C}[\lambda]^{n\times p}.
\end{align*}
With the above representation, the goal is to study how (simple) zeros of $P(\lambda)$ change to first order in $\epsilon$. Recall that the perturbed rational matrix \eqref{eq:ptfm} is the transfer function matrix of the perturbed polynomial system matrix \eqref{eq_polsysmat_per} and, by the discussion in Section \ref{sec:pert}, it is clear that it is equivalent to define such a structured condition number for the PEP referring to simple eigenvalues of (structured perturbations of) $P(\la)$ or simple zeros of the corresponding perturbations of its transfer function $R(\la)$.

\begin{definition}[Absolute structured condition number of a simple zero] \label{def:condnumber}
 Let $R(\la)$ be a regular rational matrix as in \eqref{eq:rational matrix}, and let $P(\la)$ be a regular polynomial system matrix for $R(\la)$ as in \eqref{eq_polsysmat}. Assume that $\la_0\in\C $ is a simple finite eigenvalue of $R(\la)$ and that $P(\la)$ is minimal at $\la_0$. Then, we define the structured absolute condition number of the eigenvalue $\la_0$ as follows:

\begin{align*}
\kappa_S(\lambda_0):= \lim_{\epsilon\rightarrow 0}\sup \left\{\dfrac{|\Delta \lambda_0|}{\epsilon} \,:\, [P(\la_0 + \Delta\la_0)+\epsilon\Delta P(\la_0+\Delta\la_0)][v+\epsilon\Delta v]=0 \quad \mathrm{with} \right.  &\\ 
\left.\phantom{\frac{|\Delta \lambda|}{\epsilon |\lambda|}}  \|\Delta D_i\|_2 \leq  d_i, \,\|\Delta A_i\|_2\leq  a_i\,, \|\Delta B_i\|_2\leq b_i\,, \|\Delta C_i\|_2\leq  c_i 
\right\}&,
\end{align*}
where the $a_i$, $b_i$, $c_i$, $d_i$ are nonnegative parameters, and $v$ is a right eigenvector of $P(\la)$ associated with $\la_0.$
Equivalently, $\kappa_S(\lambda_0)$ can be defined by using the transfer function matrix $R(\la)$ of $P(\la)$ as follows:
\begin{align*}
\kappa_S(\lambda_0):= \lim_{\epsilon\rightarrow 0}\sup \left\{\dfrac{|\Delta \lambda_0|}{\epsilon} \,:\, \la_0 + \Delta \la_0 \ \mathrm{is} \ \mathrm{a} \ \mathrm{zero} \ \mathrm{of }\, R(\la)+\Delta R(\la,\epsilon) \quad \mathrm{with} \right.  &\\ 
\left.\phantom{\frac{|\Delta \lambda|}{\epsilon |\lambda|}}  \|\Delta D_i\|_2 \leq  d_i, \,\|\Delta A_i\|_2\leq  a_i\,, \|\Delta B_i\|_2\leq b_i\,, \|\Delta C_i\|_2\leq  c_i 
\right\}&.
\end{align*}

\end{definition} 

The nonnegative parameters $a_i,b_i,c_i,d_i$ will be called weights and provide some freedom in how perturbations to the rational matrix $R(\lambda)$ are measured. For example, two natural choices are either $a_i=b_i=c_i=d_i=1$ (the coefficients are allowed to be perturbed by perturbations of uniform absolute magnitude) or $a_i=\| A_i \|,b_i=\|B_i\|,c_i=\|C_i\|,d_i=\|D_i\|$ (the coefficients are allowed to be perturbed by perturbations of uniform relative magnitude). We prefer to keep the weights as free parameters in our analysis as, in principle, a particular problem may suggest more peculiar choices.  In Theorem \ref{thm:condition numberlimit} we derive a formula for $\kappa_S(\lambda_0)$.

\begin{remark}\rm \label{rem:structuredcond}For any choice of positive weights, Definition \ref{def:condnumber} yields an absolute structured condition number \eqref{eq:abscond} for the eigenvalue-computing function $f : P(\la) \mapsto \la_0$ defined in Section \ref{sub:conddef}. Here, $\UU$ is the vector space of polynomial system matrices of the form \eqref{eq:psm} such that $A(\la),B(\la),C(\la)$ and $D(\la)$ are matrix polynomials of degree at most, resp., $k_A,k_B,k_C $ and $k_D$. In particular, Definition \ref{def:condnumber} is the condition number induced by the choice of the norm
\begin{equation}\label{eq:normonlyforp}
  P(\la) \mapsto \| P(\la) \|=\max \left\{ \max_{i=0}^{k_A} \frac{\| A_i \|}{a_i} ,\max_{i=0}^{k_B} \frac{\| B_i \|}{b_i},\max_{i=0}^{k_C} \frac{\| C_i \|}{c_i},\max_{i=0}^{k_D} \frac{\| D_i \|}{d_i}\right\},
\end{equation}
with the convention that if, e.g., $B(\la)=0$ then $k_B=-\infty$ and $\max_{i=0}^{k_B}\frac{ \|B_i\|}{b_i}=0$. We leave as an exercise to the reader the straightforward verification that \eqref{eq:normonlyforp} is indeed a norm, for any choice of the parameters. 
\end{remark}

\begin{remark}\rm 
	A subtler problem is to assess whether Definition \ref{def:condnumber} yields a condition number for the eigenvalue-computing function $f :R(\la) \mapsto \la_0$: indeed, \eqref{eq:normonlyforp} is clearly \emph{not} a norm on the vector space of rational functions. On the other hand, it is worth observing that the set of allowed perturbation is highly structured with respect to the original REP, i.e., not every rational small perturbation of $R(\la)$ can be written in this form (in fact, in a precise geometric sense, only very few can). Indeed, note that the tangent space at $\epsilon=0$ to the set of perturbations that we consider is 
	$$ \{ \Delta D(\la) + \Delta C(\la) A(\la)^{-1} B(\la) + C(\la) A(\la)^{-1} \Delta B(\la) - C(\la) A(\la)^{-1} \Delta A(\la) A(\la)^{-1} B(\la \}$$
	which is a finite dimensional vector space over $\C$. Indeed, even neglecting the issue of the potential non-unicity of the representation, its dimension is readily verified to be $\leq pm (\eta_A+\eta_B+\eta_C+\eta_D+4)$ where $\eta_D=k_D$ if $D(\la)\neq 0$ or $\eta_D=-1$ otherwise, and similarly for $k_A,k_B,k_C$. In contrast, the vector space of $p \times m$ rational matrices has (uncountably) infinite dimension over $\C$. Nevertheless, we believe that our definition provide a sensible practical tool and that restricting to this set of perturbation is justified when working with the representation of rational matrices that we are considering in this paper.
\end{remark}

Now, we derive a formula for the condition number introduced in Definition \ref{def:condnumber}.

\begin{theorem}[Condition number formula]\label{thm:condition numberlimit}
Let $P(\la) \in \C[\la]^{(n+p)\times(n+p)}$ 
be a regular polynomial system matrix as in \eqref{eq_polsysmat}, with $A(\la)\in\mathbb{C}[\lambda]^{n\times n}$ regular, and transfer function matrix $R(\la)$. Assume that $\la_0\in\C $ is a simple eigenvalue of $R(\la)$ and that $P(\la)$ is minimal at $\la_0$. Then the structured absolute condition number $\kappa_S(\lambda_0)$ in Definition \ref{def:condnumber} is given by
	\begin{equation*}
	\kappa_S(\lambda_0) =\frac{1}{K}
	\begin{bmatrix}
	 \|y(\la_0)^TC(\la_0)A(\la_0)^{-1}\|_2 &\|y(\la_0)\|_2
	\end{bmatrix}
	S(\lambda_0)
	\begin{bmatrix}
	\|A(\la_0)^{-1}B(\la_0)x(\la_0)\|_2 \\
	\|x(\la_0)\|_2
	\end{bmatrix},
	\end{equation*}
	by using the notation in Remark \ref{rem:notation}, where $x(\la)$ and $y(\la)$ are, respectively, right and left root vectors of $R(\lambda)$ associated with $\lambda_0$, $$K:=\displaystyle\lim_{\la\rightarrow \lambda_0} |y(\la)^T R^\prime(\la)x(\la)|\quad\text{  and  }\quad S(\lambda):=\begin{bmatrix}
	\sum_{i=0}^{k_A}a_i|\lambda|^i & \sum_{i=0}^{k_B}b_i|\lambda|^i \\
	\sum_{i=0}^{k_C}c_i|\lambda|^i & \sum_{i=0}^{k_D}d_i|\lambda|^i
	\end{bmatrix}.$$
	
\end{theorem}

\begin{proof} By using Lemma \ref{lem_der} and Proposition \ref{prop:eigenvectorslimit}, and ignoring $o(\epsilon)$ terms, we get
{\small\begin{align*}
|\Delta \lambda_0|\leq & \frac{\epsilon}{K}  \left(\|x(\la_0)\|_2\|y(\la_0)\|_2\sum_{i=0}^{k_D}d_i|\lambda_0|^i\right.+\|A(\la_0)^{-1}B(\la_0)x(\la_0)\|_2\|y(\la_0)\|_2\sum_{i=0}^{k_C}c_i|\lambda_0|^i + \\
& \|x(\la_0)\|_2\|y(\la_0)^TC(\la_0)A(\la_0)^{-1} \|_2 \sum_{i=0}^{k_B}b_i|\lambda_0|^i +\\&
\left.  \|y(\la_0)^TC(\la_0)A(\la_0)^{-1}\|_2\|A(\la_0)^{-1}B(\la_0) x(\la_0)\|_2 \sum_{i=0}^{k_A}a_i|\lambda_0|^i \right)=\\ 
& \frac{\epsilon}{K} 
\begin{bmatrix}
\|y(\la_0)^TC(\la_0)A(\la_0)^{-1}\|_2 & \|y(\la_0)\|_2
\end{bmatrix}
S(\lambda_0)
\begin{bmatrix}
\|A(\la_0)^{-1}B(\la_0)x(\la_0)\|_2 \\
\|x(\la_0)\|_2
\end{bmatrix}.
\end{align*}}

This provides an upper bound for the condition number $\kappa_S(\lambda_0)$.
We need to show that this upper bound is sharp, i.e., there is a perturbation of the form that we consider and that achieves the bound.
Let us define 
$$ \mu = \begin{cases} 0 \ &\mathrm{if} \ \la_0=0;\\
\frac{\overline{\la_0}}{|\la_0|} \ &\mathrm{otherwise},
\end{cases}$$
and let us consider  the following rank-1 matrix perturbations:
\begin{itemize}
	\item $$\Delta D_i =  d_i 
	\mu^i \dfrac{y(\la_0)x(\la_0)^T}{\|x(\la_0)\|_2\|y(\la_0)\|_2},$$ for $i=0,1,\hdots,k_D.$
	\item $$\Delta C_i =  c_i 
	\mu^i \dfrac{y(\la_0) (A(\la_0)^{-1}B(\la_0)x(\la_0))^T}{\|y(\la_0)\|_2\|A(\la_0)^{-1}B(\la_0)x(\la_0)\|_2},$$ if $A(\la_0)^{-1}B(\lambda_0)x(\la_0)\neq 0$ or $\Delta C_i =0$ otherwise, for $i=0,1,\hdots,k_C.$ 
	\item $$\Delta B_i =  b_i 
	\mu^i \dfrac{(y(\la_0)^TC(\la_0)A(\la_0)^{-1})^Tx(\la_0)^T}{\|x(\la_0)\|_2\|y(\la_0)^TC(\la_0)A(\la_0)^{-1}\|_2},$$ if $y(\la_0)^T C(\la_0)A(\la_0)^{-1}\neq 0$ or $\Delta B_i =0$ otherwise, for $i=0,1,\hdots,k_B,$ and
	\item $$\Delta A_i=  -a_i 
	\mu^i \dfrac{(y(\la_0)^T C(\la_0)A(\la_0)^{-1})^T(A(\la_0)^{-1}B(\la_0)x(\la_0))^T}{\|y(\la_0)^T C(\la_0)A(\la_0)^{-1}\|_2\|A(\la_0)^{-1}B(\la_0)x(\la_0)\|_2},$$ if $A(\la_0)^{-1}B(\la_0)x(\la_0)\neq 0$ and $ y(\la_0)^T C(\la_0)A(\la_0)^{-1}\neq 0$ or $\Delta A_i =0$ otherwise, for $i=0,1,\hdots,k_A.$
\end{itemize}

It is clear that these matrix perturbations satisfy $\|\Delta D_i\|\leq  d_i$, $\|\Delta B_i\|\leq  b_i$, $\|\Delta C_i\|\leq  c_i$ and $\|\Delta A_i\|\leq  a_i$. Then we have
{\small\begin{align*}
&y(\la_0)^T\left(\epsilon \sum_{i=0}^{k_D} \lambda_0^i \Delta D_i \right)x(\la_0) =\epsilon \|x(\la_0)\|_2\|y(\la_0)\|_2  \left( \sum_{i=0}^{k_D} |\lambda_0|^i d_i \right),\\
&y(\la_0)^T\left(\epsilon\sum_{i=0}^{k_C} \lambda_0^i \Delta C_i \right)(A(\la_0)^{-1}B(\la_0)x(\la_0)) =\epsilon  \|A(\la_0)^{-1}B(\la_0)x(\la_0)\|_2\|y(\la_0)\|_2  \sum_{i=0}^{k_C}c_i|\lambda_0|^i,\\
&(y(\la_0)^T C(\la_0)A(\la_0)^{-1})\left(\epsilon\sum_{i=0}^{k_B} \lambda_0^i \Delta B_i \right)x(\la_0) =\epsilon \|x(\la_0)\|_2 \|y(\la_0)^T C(\la_0)A(\la_0)^{-1} \|_2  \sum_{i=0}^{k_B}b_i|\lambda_0|^i,
\end{align*}}
and
{\small\begin{align*}
&(y(\la_0)^T C(\la_0)A(\la_0)^{-1}) \left(\epsilon\sum_{i=0}^{k_A} \lambda_0^i \Delta A_i \right) (A(\la_0)^{-1} B(\la_0)x(\la_0) )
= \\
& \hspace{4cm} -\epsilon \|y(\la_0)^T C(\la_0)A(\la_0)^{-1}\|_2\|A(\la_0)^{-1}B(\la_0) x(\la_0)\|_2   \sum_{i=0}^{k_A}a_i|\lambda_0|^i.
\end{align*}}%
Hence, for these perturbations all the norm inequalities are satisfied as equalities and, so, the upper bound is tight.
\end{proof}

%

\subsection{Condition number formula for zeros of $R(\la)$ that are not poles}

Theorem \ref{thm:condition numberlimit} can be simplified for those zeros $\la_0$ such that $\det A(\la_0)\neq 0$. Note that, if $\det A(\la_0)\neq 0$, the polynomial system matrix $P(\la)$ is minimal at $\la_0$. It is worth emphasizing that if $\det A(\la_0)=0$ not necessarily $\la_0$ is a pole of $A(\la)$, if we do not assume minimality on $P(\la)$ at $\la_0$, but all the poles of $R(\la)$ satisfy that $\det A(\la_0)=0$. So if $\det A(\la_0)\neq 0$ then $\la_0$ is not a pole of $R(\la)$, and the root vectors in Theorem \ref{thm:condition numberlimit} can be considered to be constant.

\begin{corollary}[Condition number formula for a zero $\la_0$ with $\det A(\la_0)\neq 0$]\label{cor:condition number}  Let $R(\la)$ be a regular rational matrix as in \eqref{eq:rational matrix}, and assume that $\la_0\in\C $ is a simple eigenvalue of $R(\la)$ with $\det A(\la_0)\neq 0$. Then the structured absolute condition number $\kappa_S(\lambda_0)$ in Definition \ref{def:condnumber} is given by
	\begin{equation*}
	\kappa_S(\lambda_0) =\frac{1}{|y^TR^\prime(\la_0)x|}
	\begin{bmatrix}
	\displaystyle\|y^TC(\lambda_0)A(\la_0)^{-1}\|_2 & \|y\|_2
	\end{bmatrix}
	S(\lambda_0)
	\begin{bmatrix}
	\|A(\la_0)^{-1}B(\lambda_0)x\|_2 \\
	\|x\|_2
	\end{bmatrix},
	\end{equation*}
	where $x$ and $y$ are, respectively, right and left eigenvectors of $R(\lambda)$ associated with $\lambda_0$ and $$S(\lambda):=\begin{bmatrix}
	\sum_{i=0}^{k_A}a_i|\lambda|^i & \sum_{i=0}^{k_B}b_i|\lambda|^i \\
	\sum_{i=0}^{k_C}c_i|\lambda|^i & \sum_{i=0}^{k_D}d_i|\lambda|^i
	\end{bmatrix}.$$
\end{corollary}

\section{Comparison with Tisseur's unstructured condition number}\label{sec:comparison}

As we explain in Remark \ref{rem:structuredcond}, Definition \ref{def:condnumber} is a structured absolute condition number for the eigenvalue computing function $f:P(\la) \mapsto \la_0$, where $P(\la)$ is a polynomial system matrix as in \eqref{eq_polsysmat} for a rational transfer function matrix $R(\la)$. That is, due to the fact that we perturb each block $A(\la),$ $B(\la)$, $C(\la)$ and $D(\la)$ of $P(\la)$ and the allowed perturbations respect, separately, the degrees of each block: the approach of studying this set of structured perturbation is natural in the context of the representation of $R(\la)$ that we assume to be given. In contrast, the condition number for matrix polynomials  defined by F. Tisseur in \cite{tisseur} considers perturbations that do not respect the degrees of  particular blocks of the corresponding polynomial matrix $P(\la)$ but the overall degree of $P(\la)$. In other words, the more general theory developed in \cite{tisseur} corresponds to studying unstructured condition numbers, as opposed to our structured condition numbers. 

 In this section, we compare our structured condition number with the unstructured condition number for $P(\la)$ obtained from \cite{tisseur}. Let us first recall Tisseur's result. For that, we expand $P(\la)$
 on the monomial basis, i.e., $P(\la)=P_d\la^d+ \cdots + P_1 \la + P_0,$
where $d=\max \{k_A,k_B,k_C,k_D\}.$ Let us now perturb each matrix coefficient $P_i$ of $P(\la)$ with a perturbation $\epsilon\Delta P_i$ such that $\| \Delta P_i \|_2 \leq  p_i,$ then the condition number formula in \cite[Theorem 5]{tisseur} for a simple eigenvalue $\la_0$ of $P(\la)$ is
\begin{equation}\label{eq_tisseurcond}
	\kappa_U(\lambda_0) =\frac{\|w\|_2\|v\|_2\left(\displaystyle \sum_{i=0}^{d}|\la_0|^i p_i\right)}{|w^TP^\prime(\la_0)v|},
\end{equation}
where $w$ and $v$ are right and left eigenvectors of $P(\la)$ associated with $\la_0,$ respectively.

\begin{remark}\rm
In \cite[Theorem 5]{tisseur} the definition of condition number, and hence its formula, is scaled by dividing by $|\la_0|$. Dividing by $|\la_0|$ is more suitable if one wants to obtain a relative condition number, as done in \cite{tisseur}. However, we prefer to present our results for absolute condition numbers as this simplifies the exposition. Note that a relative condition number can always be obtained from an absolute one via \eqref{eq:relcond}. We do not divide by $|\la_0|$ in \eqref{eq_tisseurcond} so to obtain an absolute condition number, for the purpose of a fair comparison with Definition \ref{def:condnumber}.
\end{remark} 

We now note that if each block polynomial matrix  of $P(\la)$ is perturbed so that $ \|\Delta D_i\|_2 \leq  d_i, \,\|\Delta A_i\|_2\leq  a_i\,, \|\Delta B_i\|_2\leq  b_i\,$ and $ \|\Delta C_i\|_2\leq c_i$ then the sharp upper bound for $\left\| \Delta P_i \right\|$ is $    \left\| \Delta P_i \right\|^2 \leq
     \max_{v^2 + w^2 = 1}
    (a_i v + b_i w )^2 + (c_i v + d_i w )^2.$
In other words, when we compare the structured and unstructured condition numbers, we should set
\begin{align} \label{eq_sharpbound}
    p_i := \left(\max_{v^2 + w^2 = 1}
    (a_i v + b_i w )^2 + (c_i v + d_i w )^2 \right)^{1/2}
\end{align}
as this yields the smallest set of unstructured perturbations that includes the set of structured perturbations. Moreover, we can see that 
\begin{align*}
    \max \{a_i,b_i,c_i,d_i\}^2 \leq \max_{v^2 + w^2 = 1}(a_i v + b_i w )^2 + (c_i v + d_i w )^2 \leq 4 \max \{a_i,b_i,c_i,d_i\}^2,
\end{align*}
where the lower and upper bounds are reached for some values of $a_i,b_i,c_i,d_i$. A special case of particular interest is when $a_i=b_i=c_i=d_i=1$; in this case, \eqref{eq_sharpbound} yields $p_i=2$. If we bound $\| \Delta P_i \|_2$ by $p_i$ as given in (\ref{eq_sharpbound}), we then expect that the unstructured absolute condition number $\kappa_U(\lambda_0) $ in \eqref{eq_tisseurcond} is larger than the structured absolute condition number $\kappa_S(\la_0)$ in Definition \ref{def:condnumber} or Theorem \ref{thm:condition numberlimit}. This indeed follows by definition of supremum, because the set of allowed perturbations for $\kappa_S$ is a subset of the set of allowed perturbations for $\kappa_U$. 

\begin{theorem}\label{th:comparisontisseur}  Let $P(\la)$ be a polynomial system matrix as in \eqref{eq_polsysmat} with transfer function matrix $R(\la)$. Assume that $\la_0\in\C$ is a simple zero of $R(\la)$ and that $P(\la)$ is minimal at $\la_0$. Assume that $\|\Delta D_i\|_2 \leq  d_i, \,\|\Delta A_i\|_2\leq  a_i\,, \|\Delta B_i\|_2\leq  b_i\,, \|\Delta C_i\|_2\leq c_i$ and let $p_i$ be defined as in $\ref{eq_sharpbound}$. Then
	$\kappa_S(\lambda_0) \leq \kappa_U(\lambda_0) .$
\end{theorem}

%


In the following example, we see that $\kappa_U$ can be unboundedly larger than $\kappa_S$. Then, in some applications restricting to structured perturbations can make a very significant difference.

\begin{example} Consider the rational matrix \begin{equation*}R(\lambda; \alpha, \beta, k):=
	\begin{bmatrix}
	\lambda - \alpha & 0 \\
	0 & 1
	\end{bmatrix}
	\begin{bmatrix}
	\lambda^k & \beta \\
	\beta & \lambda^k
	\end{bmatrix}^{-1}\end{equation*}
	for some parameters $\alpha, \beta$ and $ k$ such that $\alpha \neq 0$ and $k > 0$. We can set $A(\lambda; \beta, k):=
	\begin{bmatrix}
	\lambda^k & \beta \\
	\beta & \lambda^k
	\end{bmatrix},$ $
	C(\lambda; \alpha):=
	\begin{bmatrix}
	\lambda - \alpha & 0 \\
	0 & 1
	\end{bmatrix},$ $B:=
	I_2$ and $
	D:=0_2$
	such that $$P(\lambda; \alpha, \beta, k):=
	\begin{bmatrix}
	-A(\lambda; \beta, k) & B \\
	C(\lambda; \alpha) & D
	\end{bmatrix}$$
	is a polynomial system matrix of $R(\lambda; \alpha, \beta, k)$. Note that $R(\lambda; \alpha, \beta, k)$ has a zero at $\alpha$, which is also a zero of $P(\lambda; \alpha, \beta, k)$. A right eigenvector for $P$ associated with $\alpha$ is $v = \begin{bmatrix}
	1 & 0 &\alpha^k & \beta
	\end{bmatrix}^T$ and an associated left eigenvector is $w = \begin{bmatrix} 0 & 0 &1 & 0 \end{bmatrix}^T$. Then, we compute the condition numbers $\kappa_U(\alpha)$ and $\kappa_S(\alpha)$ given in \eqref{eq_tisseurcond} and Theorem \ref{thm:condition numberlimit}, respectively. We remark that, although the expression for $\kappa_S$ in Theorem \ref{thm:condition numberlimit} involves the inverse of $A$ and right and left root vectors of $R$, $\kappa_S$ can be computed by using right and left eigenvectors of $P$ taking into account Proposition \ref{prop:eigenvectorslimit} and \eqref{eq_der}. For $\kappa_S(\alpha)$, we choose the weights $d_0=0$ and $a_i=b_i=c_i=1$ for all $i \leq k_A,k_B,k_C$, respectively, and we obtain $\kappa_S(\alpha)=|\alpha|+1$. For $\kappa_U(\alpha)$, we choose the weights $p_i$ according to (\ref{eq_sharpbound}) for all $i \leq k$. This yields $p_0 = (1+ \sqrt{5})/2$, $p_1 = \sqrt{2}$ and $p_i = 1$ for $1 < i \leq k $ if $k>1$. We obtain $\kappa_U(\alpha)=\sum_{i=0}^k p_i |\alpha|^i \sqrt{1+|\alpha|^{2k}+|\beta|^2}$. Therefore, when $|\alpha| > 1$, $\kappa_U(\alpha)\rightarrow \infty$ as $k\rightarrow \infty$ or $\beta\rightarrow \infty$, while $\kappa_S(\alpha)$ keeps constant. Moreover, $\kappa_U(\alpha) / \kappa_S(\alpha) \rightarrow \infty$ as $\alpha \rightarrow \infty$. Finally, note that $\alpha$ is both a zero and a pole when $\beta = \pm \alpha^k$.
\end{example}

\section{Numerical experiments}\label{sec:num_exp}

In this section, we explore the relationship between the unstructured and structured condition numbers $\kappa_U$ and $\kappa_S$, respectively, via numerical experiments. We will show that, also in actual applications, $\kappa_U$ may be much larger than $\kappa_S$. Moreover, the choice of the representation of the rational matrix can influence the ratio between $\kappa_U$ and $\kappa_S$.

\textbf{Experiment 1.} We first consider a REP from the real-life application \eqref{eq:real_application}. One natural choice for the weights for $\kappa_S$ is to choose $a_i=b_i=c_i=d_i=1$ for all $i \leq k_A,k_B,k_C,k_D$, respectively, and otherwise the weights are set to 0. We choose the weights $p_i$ for $\kappa_U$ according to (\ref{eq_sharpbound}) for all $i \leq d$. More realistically, we might want to perturb only the actual data. In the representation given in (\ref{eq:real_application_representation}), this corresponds to setting $c_0=c_1=b_0=a_1=0$ and $d_0=d_1=b_1=a_0=1$. If we use the representation in (\ref{eq:real_application_representation}), the ratio of the condition numbers can be bounded above by 2 when using the former way of choosing the weights. In our numerical experiments, the ratio is bounded by a relatively low number also in the case when only the data is perturbed. 

\begin{figure}[h]
\centering
\begin{subfigure}{.5\textwidth}
  \centering
  \includegraphics[width=\linewidth]{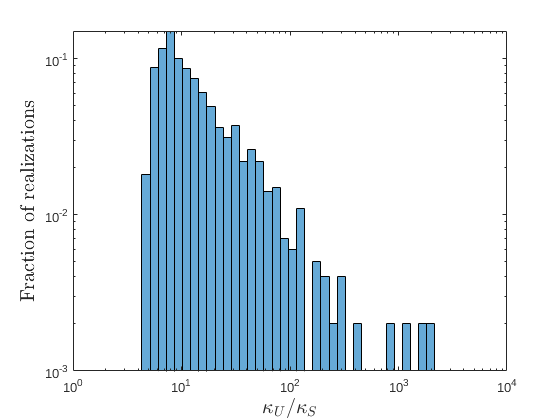}
\caption{Perturbing only actual data. }
\end{subfigure}%
\begin{subfigure}{.5\textwidth}
  \centering
  \includegraphics[width=\linewidth]{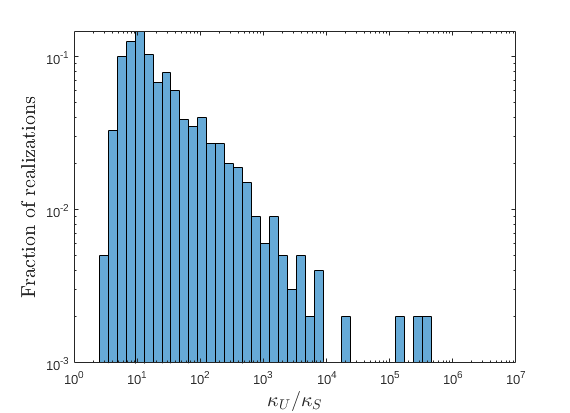}
  \caption{Perturbing all coefficient matrices. }
\end{subfigure}
\caption{For \eqref{eq:real_application_representation3}, the ratio $\kappa_U/\kappa_S$ of the condition numbers for $10^3$ realizations. }
\label{fig_real_application}
\end{figure}

It is also possible to represent (\ref{eq:real_application}) in different ways. For example, we can move one factor $\lambda$ from $B(\la)$ to $C(\la)$, or, in other words, consider the representation
\begin{equation}\label{eq:real_application_representation3}
 R(\la)=-K+\lambda M + \begin{bmatrix} \la^2 I \cdots \la^2 I
\end{bmatrix}\begin{bmatrix}
(w_1-\la) I & & \\
& \ddots & \\
& & (w_k-\la) I
\end{bmatrix}^{-1}\begin{bmatrix}
C_1 \\
\vdots \\ C_k
\end{bmatrix}
\end{equation}
rather than \eqref{eq:real_application_representation}. The structured condition number $\kappa_S$ is unchanged by this change of representation. However, the corresponding unstructured condition number $\kappa_U$ is not. In the representation \eqref{eq:real_application_representation3}, perturbing only the real data correpsonds to setting $c_0=c_1=c_2=a_1=0$ and $d_0=d_1=b_0=a_0=1$. This yields $p_0 = (1+ \sqrt{5})/2$, $p_1 = 1$, and $p_2 = 0$. Choosing $a_i=b_i=c_i=d_i=1$ for all $i \leq k_A,k_B,k_C,k_D$ yields $p_0 = 2$, $p_1 = (1+ \sqrt{5})/2$, and $p_2 = 1$. For this representation, we compute the condition numbers $\kappa_U$ and $\kappa_S$ given in (\ref{eq_tisseurcond}) and in Theorem \ref{thm:condition numberlimit}, respectively, for the largest finite eigenvalue in magnitude. We set $K, M, C_i \in \C^{200 \times 200}$ and $k=2$, where we generate the entries of the matrices as well as the poles $\omega_i$ randomly from a normal distribution of mean 0 and variance 1. We generate $10^3$ realizations in the described way and compute the ratio for the condition numbers. The results are shown in Figure \ref{fig_real_application}. We can see that, in both cases, the ratio can grow large although higher ratios are increasingly less probable.

\begin{figure}[h]
	\centering
	\begin{subfigure}{.5\textwidth}
		\centering
		\includegraphics[width=\linewidth]{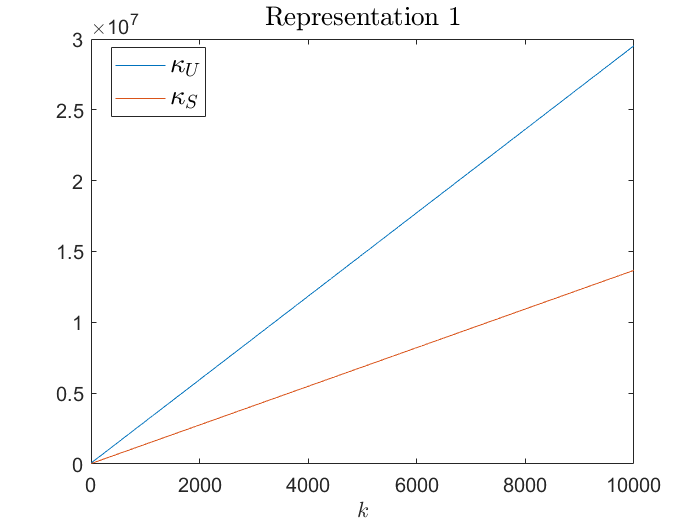}
	\end{subfigure}%
	\begin{subfigure}{.5\textwidth}
		\centering
		\includegraphics[width=\linewidth]{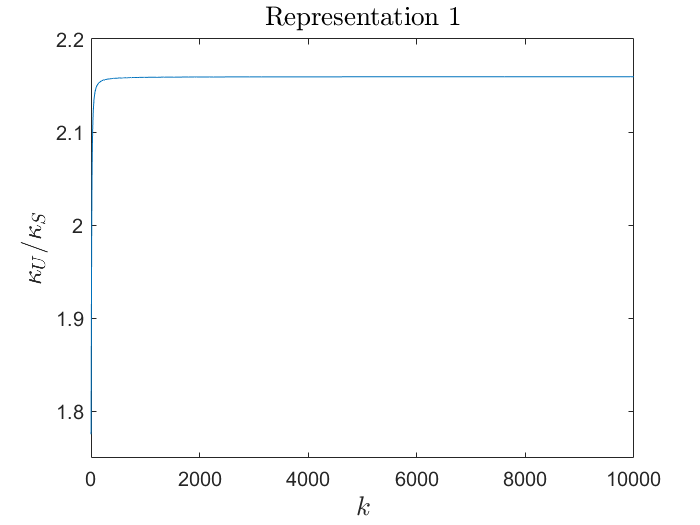}
	\end{subfigure}
	\caption{For \eqref{eq:real_application_representation2}, the condition numbers $\kappa_U$ and $\kappa_S$  as well as their ratio as a function of the parameter $k$.  }
	\label{fig_loaded_string2}
\end{figure}

\begin{figure}[h]
\centering
\begin{subfigure}{.5\textwidth}
  \centering
  \includegraphics[width=\linewidth]{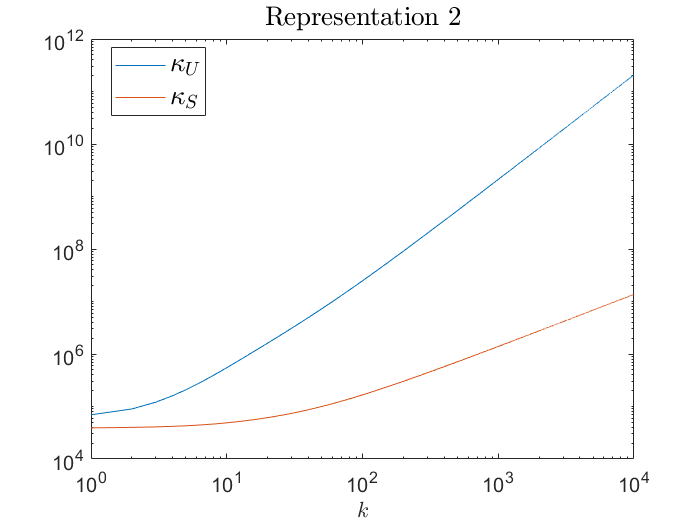}
\end{subfigure}%
\begin{subfigure}{.5\textwidth}
  \centering
  \includegraphics[width=\linewidth]{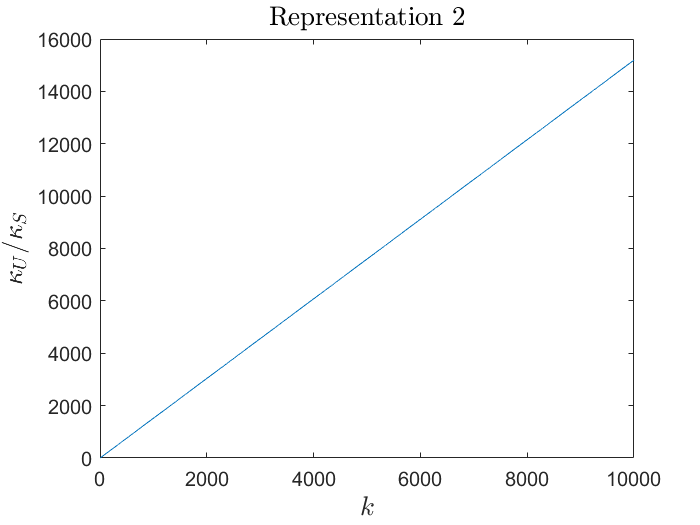}
\end{subfigure}
\caption{For \eqref{eq:real_application_representation4}, the condition numbers $\kappa_U$ and $\kappa_S$  as well as their ratio as a function of the parameter $k$. }
\label{fig_loaded_string1}
\end{figure}

\textbf{Experiment 2.} Let us next consider the real-life application \eqref{eq:real_application2}. In particular, we consider two different ways to represent this rational matrix. The first representation is \eqref{eq:real_application_representation2} while the second representation is \begin{equation}\label{eq:real_application_representation4}
R(\la)=A-\la B + k\la e_n\begin{bmatrix}
\la-k/m
\end{bmatrix}^{-1}e_n^T.
\end{equation} 
 For perturbations, we set $a_i=b_i=c_i=d_i=1$ for all $i \leq k_A,k_B,k_C,k_D$, respectively, and otherwise we set the weights to 0. Again, we choose the weights $p_i$ for $\kappa_U$ according to (\ref{eq_sharpbound}). Let us fix $n=10$ and $m=1$, and let $k$ vary. We compute the condition numbers $\kappa_U$ and $\kappa_S$ for the largest finite eigenvalue in magnitude for both representations. The results are in Figures \ref{fig_loaded_string2} and \ref{fig_loaded_string1}.

In Figure \ref{fig_loaded_string2}, both $\kappa_U$ and $\kappa_S$ are linear in $k$. However, in Figure \ref{fig_loaded_string1} the slope for $\kappa_U$ in the loglog-plot is close to 2, which implies that $\kappa_U$ is quadratic in $k$; whereas the slope for $\kappa_S$ is close to 1, which implies that $\kappa_S$ is linear in $k$. As a result, the ratio $\kappa_U / \kappa_S$ diverges as $k$ approaches infinity.
In addition, we have that the structured condition number $\kappa_S$ is the same for the two representations \eqref{eq:real_application_representation2} and \eqref{eq:real_application_representation4}, whereas the corresponding unstructured condition number $\kappa_U$ does depend on this choice. 

\section{Conclusions and open problems}\label{sec:con}
	
	We defined a structured condition number $\kappa_S$ for eigenvalues $\la_0$ of a (locally) minimal polynomial system matrix, which are the zeros of the corresponding transfer function $R(\la)$. We then derived a computable expression for $\kappa_S$, that is valid even in the case of $\la_0$ being also a pole of $R(\la)$. For that, we used the notion of root vectors. We finally compared our structured condition number $\kappa_S$ with the standard condition number for polynomial matrices defined by F. Tisseur in \cite{tisseur}, showing that $\kappa_S$ is never larger and that may actually be much smaller in some applications, taking into account different representations of $R(\la)$. 
	
	An open problem is to study which representations of rational matrices have favourable properties with respect to conditioning. In addition, nowadays the standard approach to solve the REP is via linearizations. That is, transforming the REP into a generalized eigenvalue problem in such a way that the spectral information is preserved. It is also an open to study which linearizations have favourable properties with respect to conditioning. Another related problem is to derive a computable expression for local backward errors for pairs of approximate eigenvalues and eigenvectors of rational matrices.

\section*{Acknowledgements}

VN and MQ thank Froil\'{a}n Dopico for enlighetning comments on the topic of this paper and for chupitos, both offered during a research visit at Universidad Carlos III de Madrid.

\end{document}